\documentclass[a4paper,11pt,english]{smfart}

\usepackage{amssymb}
\usepackage{url}
\usepackage[T1]{fontenc}
\RequirePackage{calrsfs}
\DeclareSymbolFont{rsfscript}{OMS}{rsfs}{m}{b}
\DeclareSymbolFontAlphabet{\mathrsfs}{rsfscript}
\usepackage[utopia,expert]{mathdesign}
\usepackage{lscape}
\usepackage{fancybox}
\usepackage{makecell}
\usepackage{todonotes}

\usepackage[leftbars]{changebar}

\usepackage{palatino}
\usepackage{rotating}
\usepackage{graphicx}

\usepackage{enumerate}

\usepackage{color}
\definecolor{shadecolor}{gray}{0.90}

\def\bfit{\bfseries\itshape}

\def\proofname{D\'emonstration}

\input xypic
\xyoption{all}
\xyoption{arc}

\makeindex

\def\indexname{Index des notations}

\newtheorem{theo}{Theorem}[section]
\def\thetheo{\thesection.\arabic{theo}}
\newtheorem{prop}[theo]{Proposition}

\newtheorem{lem}[theo]{Lemma}

\newtheorem{coro}[theo]{Corollary}

\renewcommand\theequation{\thetheo}
\def\equat{\refstepcounter{theo}\begin{equation}}
\def\endequat{\end{equation}}

\renewcommand\thetable{\thetheo}

\renewcommand\thesection{\arabic{section}}
\renewcommand\thesubsection{\thesection.\Alph{subsection}}

\def\contentsname{Table des mati\`eres}
\def\listfigurename{Liste des figures}
\def\listtablename{Liste des tables}
\def\appendixname{Appendice}

\def\refname{R\'ef\'erences}

    \def\CM{{\mathbb{C}}}
    
  \def\eG{{\mathfrak e}}  
    
  \def\gG{{\mathfrak g}}  
    
  \def\iG{{\mathfrak i}}

\def\LG{{\mathfrak L}}  \def\lG{{\mathfrak l}}  
  \def\mG{{\mathfrak m}}

\def\QG{{\mathfrak Q}}  \def\qG{{\mathfrak q}}  
    
  \def\sG{{\mathfrak s}}

\def\ZG{{\mathfrak Z}}    \def\ZM{{\mathbb{Z}}}

  \def\ab{{\mathbf a}}

  \def\eb{{\mathbf e}}  \def\EC{{\mathcal{E}}}
    
\def\Gb{{\mathbf G}}    
\def\Hb{{\mathbf H}}

\def\Lb{{\mathbf L}}    \def\LC{{\mathcal{L}}}
    \def\MC{{\mathcal{M}}}

\def\Sb{{\mathbf S}}    
    
  \def\ub{{\mathbf u}}

    \def\ZC{{\mathcal{Z}}}

  \def\mrm{{\mathrm{m}}}

\def\Srm{{\mathrm{S}}}

  \def\yrm{{\mathrm{y}}}  
\def\Zrm{{\mathrm{Z}}}

          \def\mdo{{\dot{m}}}

\def\a{\alpha}
\def\b{\beta}
\def\g{\gamma}
\def\G{\Gamma}
\def\d{\delta}

\def\e{\varepsilon}
\def\ph{\varphi}

\def\L{\Lambda}

\def\O{\Omega}
\def\r{\rho}
\def\s{\sigma}

\def\th{\theta}
\def\Th{\Theta}
\def\t{\tau}

\def\z{\zeta}

\def\delb{{\boldsymbol{\delta}}}

\def\mub{{\boldsymbol{\mu}}}

\DeclareMathOperator{\diag}{{\mathrm{diag}}}

\DeclareMathOperator{\Ker}{{\mathrm{Ker}}}

\def\to{\rightarrow}
\def\longto{\longrightarrow}

\def\fonction#1#2#3#4#5{\begin{array}{rccc}
{#1} : & {#2} & \longto & {#3}  \\
& {#4} & \longmapsto & {#5} 
\end{array}}

\def\fonctio#1#2#3#4{\begin{array}{ccc}
{#1} & \longto & {#2} \\
{#3} & \longmapsto & {#4} 
\end{array}}

\def\DS{\displaystyle}
\def\SS{\scriptstyle}

\def\finl{~$\blacksquare$}

\def\lexp#1#2{\kern\scriptspace\vphantom{#2}^{#1}\kern-\scriptspace#2}
\def\le{\hspace{0.1em}\mathop{\leqslant}\nolimits\hspace{0.1em}}
\def\ge{\hspace{0.1em}\mathop{\geqslant}\nolimits\hspace{0.1em}}

\mathchardef\inferieur="321E
\mathchardef\superieur="321F

\def\eqna{\begin{eqnarray*}}
\def\endeqna{\end{eqnarray*}}

\catcode`\@=11
\long\def\@car#1#2\@nil{#1}
\long\def\@first#1#2{#1}
\long\def\@second#1#2{#2}
\long\def\ifempty#1{\expandafter\ifx\@car#1@\@nil @\@empty
  \expandafter\@first\else\expandafter\@second\fi}
\catcode`\@=12

\def\Ref{{\mathrm{Ref}}}

\def\boitegrise#1#2{\begin{centerline}{\fcolorbox{black}{shadecolor}{~
    \begin{minipage}[t]{#2}{\vphantom{~}#1\vphantom{$A_{\DS{A_A}}$}}
            \end{minipage}~}}\end{centerline}\medskip}

\theoremstyle{remark}
\newtheorem{rema}[theo]{Remark}

\theoremstyle{plain}

\def\BIL{LR}
\def\GAUCHE{L}
\def\CAR{CAR}
\def\FAM{FAM}

\def\euler{{\eb\ub}}

\def\xyinj{\ar@{^{(}->}}
\def\xysur{\ar@{->>}}

\bigskip

\makeatletter
\def\hlinewd#1{%
\noalign{\ifnum0=`}\fi\hrule \@height #1 %
\futurelet\reserved@a\@xhline}
\makeatother

\newlength\epaisLigne

\usepackage{multirow}

\newcommand{\longiso}{\stackrel{\sim}{\longrightarrow}}

\makeatletter
\def\hlinewd#1{%
\noalign{\ifnum0=`}\fi\hrule \@height #1 %
\futurelet\reserved@a\@xhline}
\makeatother

\usepackage{multirow}

\usepackage{lscape}

\def\troncation{{\mathrm{Trunc}}}

\def\itembul{\item[$\bullet$]}

\def\Lie{{{\LG\iG\eG}}}

\begin{document}


\title{On the Calogero-Moser space associated \\
with dihedral groups II. \\
The equal parameter case}

\author{{\sc C\'edric Bonnaf\'e}}
\address{IMAG, Universit\'e de Montpellier, CNRS, Montpellier, France} 

\makeatletter
\email{cedric.bonnafe@umontpellier.fr}
\makeatother

\date{\today}

\thanks{The author is partly supported by the ANR:
Projects No ANR-16-CE40-0010-01 (GeRepMod) and ANR-18-CE40-0024-02 (CATORE)}

\maketitle

\pagestyle{myheadings}

\markboth{\sc C. Bonnaf\'e}{\sc Calogero-Moser space and dihedral groups: 
the equal parameter case}


\begin{abstract}
We continue the study of Calogero-Moser spaces associated 
with dihedral groups by investigating in more details the 
equal parameter case: we obtain explicit equations, 
some informations about the Poisson bracket, the structure 
of the Lie algebra associated with the cuspidal point and 
the action of $\Sb\Lb_2(\CM)$.
\end{abstract}

\bigskip

We continue here the study of Calogero-Moser space $\ZC_c$ associated 
with the dihedral group $W$ of order $2d$ started in~\cite{bonnafe diedral}, from 
which we keep the notation. We mainly focus on the equal 
parameter case (i.e. the case where $a=b$ with the notation 
of~\cite[\S{3.4}]{bonnafe diedral})\footnote{Recall that, if $d$ is odd, 
then we have necessarily $a=b$.}. In this case, the main results of this 
paper are the following:
\begin{itemize}
\itembul We describe explicit equations for $\ZC_c$.

\itembul We obtain informations about the Poisson bracket 
that allow to determine the structure 
of the Lie algebra associated with the cuspidal point.

\itembul We describe the action of $\Sb\Lb_2(\CM)$ on 
the generators of $Z_c$ and explain how the presentation of $Z_0$ 
can be interpreted in terms of Hermite's reciprocity 
law\footnote{We wish to thank warmly Pierre-Louis Montagard 
for his enlighting explanations.} (see for instance~\cite[Cor.~2.2]{brion}).

\itembul If $\t$ denotes the diagram automorphism of $W$, then $\t$ acts on 
$\ZC_c$ because we are in the equal parameter case, and we prove that 
the irreducible components of $\ZC_c^\t$ are also Calogero-Moser spaces 
associated with other reflection groups. This confirms~\cite[Conj.~FIX]{calogero} 
(or~\cite[Conj.~B]{bonnafe auto}) in this small case.
\end{itemize}
These results will be used by 
G. Bellamy, B. Fu, D. Juteau, P. Levy, E. Sommers and the author in a forthcoming 
paper, where it will be shown that, for $d \ge 5$, the 
symplectic singularity of $\ZC_c$ at its cuspidal point 
is a new family of isolated symplectic singularities 
whose local fundamental group is trivial~\cite{bbjfls symplectic}, 
answering an old question of Beauville~\cite{beauville}.

\bigskip

These computations are based on a first paper of the author 
on Calogero-Moser spaces associated with dihedral 
groups~\cite{bonnafe diedral} and on an algorithm developed 
by U. Thiel and the author~\cite{bonnafe thiel}. This algorithm 
was implemented by Thiel~\cite{thiel} 
in his {\sc Champ} package for {\sc Magma}~\cite{magma}.
Explicit computer computations in small cases (i.e. $d \in \{4,5,6,7\}$) 
were necessary to find the general pattern. So, even though 
this does not appear in this paper, 
it is fair to say that the above results owe their existence to {\sc Magma}.

\bigskip

\noindent{\bfit Recollection of notation from~\cite{bonnafe diedral}.---} 
We will use the notation of the first part~\cite{bonnafe diedral} and we 
recall here some of them, the most important ones. 
We set $V=\CM^2$ and $(x,y)$ denotes its canonical basis while 
$(X,Y)$ is the dual basis of $V^*$. We identify $\Gb\Lb_\CM(V)$ with $\Gb\Lb_2(\CM)$. 
We also fix a non-zero natural number $d$, as well as a primitive 
$d$-th root of unity $\z \in \CM^\times$. If $i \in \ZM$ or $\ZM/d\ZM$, we set
$$s_i=\begin{pmatrix} 0 & \z^i \\ \z^{-i} & 0 \end{pmatrix},$$
$s=s_0$, $t=s_1$ and $W=\langle s,t \rangle$: it is the dihedral group 
of order $2d$. The set $\Ref(W)$ of reflections of $W$ is $\{s_i~|~i \in \ZM/d\ZM\}$. 
Finally, let $w_0$ denote the longest element of $W$ (we have $w_0=t(st)^{(d-1)/2}$ 
if $d$ is odd and $w_0=(st)^{d/2}$ if $d$ is even): this notation was used 
in the first part~\cite[Rem.~6.4]{bonnafe diedral} but we had forgotten 
to define it! It will be used here in Section~\ref{sec:fixe}.

We set $q=xy$, $Q=XY$, $r=x^d+y^d$, $R=X^d+Y^d$ and, if $0 \le i \le d$, 
$$\ab_{i,0} = x^{d-i} Y^i + y^{d-i}Y^i.$$
In this second part, we will not use the notation $r$ or $R$ as $r=\ab_{0,0}$ 
and $R=\ab_{d,0}$: we prefer this second notation. If $i \ge 0$, we set 
$$\euler_0^{(i)}=(xX)^i+(yY)^i$$
and $\euler_0=\euler_0^{(1)}$. 

We fix a map $c : \Ref(W) \to \CM$ and we set $a=c_s$ and $b=c_t$. We denote by 
$\Hb_c$ the rational Cherednik algebra at $t=0$, with parameter $c$, whose 
presentation is given in~\cite[$($3.2$)$]{bonnafe diedral}. Its center is denoted by $Z_c$ 
and we denote by $\ZC_c$ the affine variety whose algebra of regular functions 
$\CM[\ZC_c]$ is precisely $Z_c$. 

We denote by $\troncation_c$ the $\CM$-linear map
$$\troncation_c : \Hb_c \longto \CM[V \times V^*]$$
such that, if $f \in \CM[V \times V^*]$ and $w \in W$, then 
$$\troncation_c(f w)=
\begin{cases}
f & \text{if $w =1$,}\\
0 & \text{otherwise.}
\end{cases}$$
It is the map induced by the map $\troncation$ 
defined in~\cite[\S{3.4}]{bonnafe diedral}. Its restriction 
$\troncation_c : Z_c \to \CM[V \times V^*]^W$ is an isomorphism 
of $\ZM$-graded vector spaces~\cite[Lem.~3.5]{bonnafe diedral}.
Recall that it is $P_{\!\bullet}$-linear, 
where $P_{\! \bullet}=\CM[V]^W \otimes \CM[V^*]^W=\CM[q,Q,\ab_{0,0},\ab_{d,0}]$.

We add a further notation which will be useful in this second part, namely, we set 
$$e=\begin{pmatrix} 0 & 1 \\ 0 & 0 \end{pmatrix},
\quad h=\begin{pmatrix} 1 & 0 \\ 0 & -1 \end{pmatrix}
\quad \text{and}\quad f=\begin{pmatrix} 0 & 0 \\ 1 & 0 \end{pmatrix},$$
so that $(e,h,f)$ is the standard basis of the Lie algebra $\sG\lG_2(\CM)$.

\bigskip

\boitegrise{{\bf Hypothesis.} 
{\it All along this paper, together with the above notation, 
we make the additional assumption that $a=b$. 
Recall that it is automatically satisfied if $d$ is odd.}}{0.75\textwidth}

\section{Back to ${\boldsymbol{\ZC_0=(V \times V^*)/W}}$} 

\medskip

\subsection{Some polynomial identities}
If $i \ge 0$, let $\euler_0^{[i]}$ denote the element
$$\euler_0^{[i]}=\frac{(xX)^{i+1}-(yY)^{i+1}}{xX-yY}=\sum_{j=0}^i (xX)^{i-j}(yY)^j$$
of $Z_0=\CM[V \times V^*]^W$. In other words, with the notation of~\cite[\S{2}]{bonnafe diedral}, 
$$\euler_0^{[i]}=\sum_{0 \le j < i/2} (qQ)^j \euler_0^{(i-2j)} + \delb_{\text{$i$ is even}},$$
where $\delb_{\text{$i$ is even}}$ is equal to $1$ (resp. $0$) 
if $i$ is even (resp. odd). Hence, using the inversion formula~\cite[(2.1)]{bonnafe diedral}, one gets 
$$\euler_0^{[i]}=\sum_{0 \le j < i/2} \Bigl((qQ)^j 
\sum_{0 \le k \le 1/2-j} n_{i-2j,k} (qQ)^k \euler_0^{i-2j-2k}\Bigr) + \delb_{\text{$i$ is even}},$$
which can be rewritten
\equat\label{eq:eu-i}
\euler_0^{[i]}=\sum_{0 \le j \le i/2} m_{i,j} (qQ)^j \euler_0^{i-2j},
\endequat
for some elements $m_{i,j} \in \ZM$. Let $\Psi_i(T,T',T'')$ denote the 
polynomial in three indeterminates equal to 
$\sum_{0 \le j \le i/2} m_{i,j} (T'T'')^j T^{i-2j}$. It is homogeneous of degree $i$ 
for the natural graduation of $\CM[T,T',T'']$ and, as a polynomial in $T$ with coefficients 
in $\CM[T',T'']$, it is monic. If we denote by $\CM[T,T',T'']_k$ the homogeneous 
component of $\CM[T,T',T'']$ of degree $k$, then~(\ref{eq:eu-i}) shows that 
\equat\label{eq:base}
\text{\it $(T^{\prime k-j} T^{\prime\prime i} \Psi_{j-i})_{0 \le i \le j \le k}$ is a basis 
of $\CM[T,T',T'']_k$.}
\endequat
By construction, $\Psi_i$ is the unique polynomial 
satisfying the following identity:
\equat\label{eq:pi}
\Psi_i(\euler_0,q,Q)=\frac{(xX)^{i+1}-(yY)^{i+1}}{xX-yY}.
\endequat
The unicity comes from the fact that $\euler_0$, $q$ and $Q$ are algebraically 
independent.
Note that $\Psi_0=1$ and $\Psi_1=T$. Now the sequence $(\Psi_i)_{i \ge 0}$ is 
easily determined by the following recursive formula: if $i \ge 1$, then 
\equat\label{eq:pi-rec}
\Psi_{i+1}=T \Psi_i - T'T'' \Psi_{i-1}.
\endequat
Indeed, this follows from the fact that 
$(xX)^{i+2}-(yY)^{i+2}=(xX+yY)((xX)^{i+1}-(yY)^{i+1})-xyXY((xX)^i-(yY)^i)$.  
Note also for future reference the following two relations: if $i \ge 1$, then 
\equat\label{eq:pi-der}
\begin{cases}
\DS{2T' \frac{\partial \Psi_i}{\partial T} + T \frac{\partial \Psi_i}{\partial T''} = (i+1)T'\Psi_{i-1}},\\
~\\
\DS{2T'' \frac{\partial \Psi_i}{\partial T} + T \frac{\partial \Psi_i}{\partial T'} = (i+1)T''\Psi_{i-1}}.
\end{cases}
\endequat

\bigskip

\begin{proof}[Proof of~\eqref{eq:pi-der}]
We prove only the first identity, the second one being obtained by exchanging 
the roles of $(x,y)$ and $(X,Y)$. 
Let us consider the two identities obtained by applying $\partial/\partial X$ and 
$\partial/\partial Y$ to~(\ref{eq:pi}):
$$
\begin{cases}
\DS{x \frac{\partial \Psi_i}{\partial T}(\euler_0,q,Q)+ Y \frac{\partial \Psi_i}{\partial T''}(\euler_0,q,Q)
= \frac{(i+1)x^{i+1}X^i(xX-yY)-x((xX)^{i+1}-(yY)^{i+1})}{(xX-yY)^2}},\\
~\\
\DS{y \frac{\partial \Psi_i}{\partial T}(\euler_0,q,Q)+ X \frac{\partial \Psi_i}{\partial T''}(\euler_0,q,Q)
= \frac{-(i+1)y^{i+1}Y^i(xX-yY)+y((xX)^{i+1}-(yY)^{i+1})}{(xX-yY)^2}.}\\
\end{cases}
$$
Multiplying the first equality by $y$, the second by $x$, and adding the results yields exactly 
$$2q \frac{\partial \Psi_i}{\partial T}(\euler_0,q,Q)+\euler_0 \frac{\partial \Psi_i}{\partial T''}(\euler_0,q,Q)
= (i+1)q\Psi_{i-1}(\euler_0,q,Q),$$
as expected.
\end{proof}

\subsection{Presentation} 
We rewrite slightly differently the presentation 
of $Z_0=\CM[V \times V^*]^W$ obtained in~\cite[Theo.~2.1]{bonnafe diedral} 
according to our needs. A straightforward computation shows that, if $1 \le i \le j \le d-1$, then 
$$\ab_{i-1,0}\ab_{j+1,0} - \ab_{i,0}\ab_{j,0} 
= (\euler_0^2-4qQ) q^{d-j-1} Q^{i-1} \euler_0^{[j-i]} .$$
Using~(\ref{eq:eu-i}), this gives
$$
\ab_{i-1,0}\ab_{j+1,0} - \ab_{i,0}\ab_{j,0} = (\euler_0^2-4qQ)
q^{d-j-1} Q^{i-1} \Psi_{j-i}(\euler_0,q,Q) 
\leqno{(\ZG_{i,j}^0)}$$
This equation can also be obtained by substracting the equation 
$(\Zrm_{i,j}^0)$ to the equation $(\Zrm_{i-1,j+1}^0)$ (with the notation 
of~\cite[\S{2}]{bonnafe diedral}). Consequently, 
the presentation given in~\cite[Theo.~2.1]{bonnafe diedral} can be rewritten 
as follows:

\bigskip

\begin{theo}\label{theo:presentation-z0}
The algebra of invariants $\CM[V \times V^*]^W$ admits the following presentation:
\begin{itemize}
\item[$\bullet$] Generators: $q$, $Q$, $\euler_0$, $\ab_{0,0}$, $\ab_{1,0}$, $\ab_{2,0}$,\dots, 
$\ab_{d,0}$.

\item[$\bullet$] Relations:
$$
\begin{cases}
\euler_0 \ab_{i,0} = q \ab_{i+1,0} + Q \ab_{i-1,0} & \text{\it for $1 \le i \le d-1$,}\\
\ab_{i-1,0}\ab_{j+1,0} - \ab_{i,0}\ab_{j,0} = 
(\euler_0^2-4qQ)q^{d-j-1} Q^{i-1} \Psi_{j-i}(\euler_0,q,Q)  & \text{\it for $1 \le i \le j \le d-1$.}\\
\end{cases}
$$
\end{itemize}
\end{theo}

\bigskip

\subsection{Poisson bracket} 
The Poisson bracket on $\CM[V \times V^*]^W$ is obtained by restriction 
of the natural one on $\CM[V \times V^*]$, which is completely 
determined by the following rules:
$$\{x,X\} = \{y,Y\}=1\qquad\text{and}\qquad 
\{x,y\}=\{X,Y\}=\{x,Y\}=\{y,X\}=0.$$
Therefore, a straightforward computation shows that 
the Poisson bracket between the generators of $\CM[V \times V^*]^W$ 
is given by:
\equat\label{eq:poisson-0}
\begin{cases}
\{q,Q\}= \euler_0,\\
\{\euler_0,q\} = -2q,\\
\{\euler_0,Q\}= 2Q,\\
\{\euler_0,\ab_{i,0}\}=(2i-d) \ab_{i,0},\\
\{q,\ab_{i,0}\} = i \ab_{i-1,0}\\
\{Q,\ab_{i,0}\}=(i-d)\ab_{i+1,0}\\
\{\ab_{i,0},\ab_{j,0}\} = j(d-i) q^{d-j}Q^i \euler_0^{(j-i-1)} - i(d-j) q^{d-j-1}Q^{i-1} 
\euler_0^{(j-i+1)},
\end{cases}
\endequat
where the last equality only holds if $0 \le i < j \le d$. 
In particular, $(Q,\euler_0,-q)$ is an $\sG\lG_2$-triple 
(for the Lie algebra structure on $\CM[V \times V^*]^W$ induced 
by the Poisson bracket). Note that 
\equat\label{eq:sl2-det}
\{Q,\euler_0^2-4qQ\}=\{q,\euler_0^2-4qQ\}=\{\euler_0,\euler_0^2-4qQ\}=0.
\endequat

\def\bull{{\SS{\,\bullet\,}}}

\subsection{Action of $\Sb\Lb_2(\CM)$}\label{sub:sl2-0}
Since $W$ is a Coxeter group, the $\CM W$-modules $V$ and $V^*$ are isomorphic. 
In our situation, the map
$$\fonction{\Phi}{V}{V^*}{\a x + \b y}{\b X + \a Y}$$
is an isomorphism of $\CM W$-modules. One then gets an action 
of $\Sb\Lb_2(\CM)$ on $V \times V^*$ as follows:
$$\begin{pmatrix} \a & \b \\ \g & \d \end{pmatrix} \cdot (u,U) = 
(\a u + \b \Phi^{-1}(U), \g \Phi(u) + \d U).$$
By construction, this action commutes with the action of $W$, 
so induces an action of $\Sb\Lb_2(\CM)$ on the $\CM$-algebras 
$\CM[V \times V^*]$, $\CM[V \times V^*] \rtimes W$ 
and $\CM[V \times V^*]^W$. This induces an action 
of the Lie algebra $\sG\lG_2(\CM)$ by derivations on $\CM[V \times V^*]$ and 
$\CM[V \times V^*]^W$. For conventional reason, if $\ph \in \CM[V \times V^*]$ and 
$\xi \in \sG\lG_2(\CM)$, we denote by $\xi \bull \ph$ the image of $\ph$ under the action of 
$-\lexp{t}{\xi}$. It is easily checked 
on the generators $x$, $y$, $X$, $Y$ of $\CM[V \times V^*]$ that 
\equat\label{eq:action-sl2}
e \bull \ph = \{Q,\ph\},\qquad h \bull \ph = \{\euler_0,\ph\}
\qquad\text{and}\qquad f \bull \ph =\{-q,\ph\}
\endequat
for all $\ph \in \CM[V \times V^*]$.

\bigskip

\def\abb{{\boldsymbol{a}}}
\def\bbb{{\boldsymbol{b}}}
\def\eulerb{{\boldsymbol{e\!\!u}}}

\section{Calogero-Moser space at equal parameters}

\medskip

\boitegrise{{\bf Notation.} {\it We denote 
by $q$, $Q$, $\eulerb$, $\abb_0$, $\abb_1$,\dots, $\abb_d$ the respective images of 
$q$, $Q$, $\euler$, $\ab_0$, $\ab_1$,\dots, $\ab_d$ in $Z_c$. }}{0.75\textwidth}

\medskip

Note the following formulas:
\equat\label{eq:x-com}
\begin{cases}
[x,X]=-a\DS{\sum_{i \in \ZM/d\ZM} s_i},\\
\\
[x,Y]=a\DS{\sum_{i \in \ZM/d\ZM} \z^{-i} s_i},\\
\\
[y,X]=a\DS{\sum_{i \in \ZM/d\ZM} \z^i s_i},\\
\\
[y,Y]=-a\DS{\sum_{i \in \ZM/d\ZM} s_i}.
\end{cases}
\endequat
Note also the following formula, which follows from~\cite[\S{3.6}]{gordon}: 
if $P \in \CM[X,Y]$, then
\equat\label{eq:x-com-p}
[x,P]=-a\sum_{i \in \ZM/d\ZM} \frac{P-\lexp{s_i}{P}}{X-\z^iY} s_i = 
-a\sum_{i \in \ZM/d\ZM} s_i \frac{P-\lexp{s_i}{P}}{X-\z^iY}.
\endequat

\subsection{Explicit form of the generators} 
The elements $\eulerb$, $\abb_0$, $\abb_1$,\dots, $\abb_d$ are characterized 
by the fact that $\troncation_c(\eulerb)=\euler_0$ and 
$\troncation_c(\abb_j)=\ab_{j,0}$.
Recall from~\cite[\S{3.3}~and~\S{4.1}]{calogero} that 
\equat\label{eq:euler}
\eulerb=xX+yY+a \sum_{i \in \ZM/d\ZM} s_i.
\endequat
An important feature of the equal parameter case is that the elements $\abb_j$ 
have a reasonably simple form:

\bigskip

\begin{prop}\label{prop:ai}
If $0 \le j \le d$, then
\eqna
\abb_j &=& \DS{x^{d-j} Y^j + y^{d-j} X^j - a \sum_{i \in \ZM/d\ZM} \z^{-ij} 
\,\frac{x^{d-j}-\z^{ij} y^{d-j}}{x-\z^{-i}y} \cdot \frac{X^j-\z^{ij} Y^j}{X-\z^iY}\, s_i} \\
&=& \DS{x^{d-j} Y^j + y^{d-j} X^j - a \sum_{i \in \ZM/d\ZM} \z^{-ij} 
\,\frac{x^{d-j}-\z^{ij} y^{d-j}}{x-\z^{-i}y} ~s_i~ \frac{X^j-\z^{ij} Y^j}{X-\z^iY}} \\
\endeqna
\end{prop}

\bigskip

\noindent{\bfit Notation.---} For future use of the above formula, we set
$$\g_{i,j}=\frac{x^{d-j}-\z^{ij} y^{d-j}}{x-\z^{-i}y}\qquad\text{and}\qquad 
\G_{i,j}=\frac{X^j-\z^{ij} Y^j}{X-\z^iY}$$
for $i \in \ZM/d\ZM$ and $0 \le j \le d$. Note that $\g_{i,d}=\G_{i,0}=0$.\finl

\bigskip

\begin{proof}
Let $\bbb_j \in \Hb_c$ denote the right-hand side of the equation 
of the proposition. 
Since $\troncation_c$ induces an isomorphism $Z_c \longiso \CM[V \times V^*]^W$ 
and $\troncation_c(\abb_j)=\troncation_c(\bbb_j)$, it is sufficient to check that $\bbb_j \in Z_c$. 
First, an easy computation shows that $\bbb_j$ commutes with $s$ and $t$. Now, 
by~(\ref{eq:x-com-p}), we have
\eqna
[x,\bbb_j]&=&
\DS{x^{d-j}\Bigl(a \sum_{i \in \ZM/d\ZM} \z^{-ij} s_i \frac{X^j-\z^{ij} Y^j}{X - \z^i Y} \Bigr)
+y^{d-j} \Bigl(-a \sum_{i \in \ZM/d\ZM} s_i \frac{X^j-\z^{ij}Y^j}{X-\z^iY}  \Bigr) }\\
&&\DS{-a \sum_{i \in \ZM/d\ZM} \z^{-ij} \frac{x^{d-j}-\z^{ij} y^{d-j}}{x-\z^{-i}y} 
\Bigl[ x,s_i \frac{X^j-\z^{ij}Y^j}{X-\z^i Y} \Bigr] } \\
&=& \DS{a\sum_{i \in \ZM/d\ZM}(\z^{-ij}x^{d-j}-y^{d-j}) s_i \frac{X^j-\z^{ij}Y^j}{X-\z^i Y} }\\
&&\DS{- a \sum_{i \in \ZM/d\ZM} \frac{\z^{-ij}x^{d-j}-y^{d-j}}{x-\z^{-i}y} [x,s_i] 
\frac{X^j-\z^{ij}Y^j}{X-\z^i Y} } \\
&&\DS{- a \sum_{i \in \ZM/d\ZM} \z^{-ij} \frac{x^{d-j}-\z^{ij} y^{d-j}}{x-\z^{-i}y} s_i 
\Bigl[x,\frac{X^j-\z^{ij}Y^j}{X-\z^i Y}\Bigr] }
\endeqna
Now, the first two lines of this last equation compensate each other and it remains
\eqna
[x,\bbb_j]&=&\DS{-a \sum_{i \in \ZM/d\ZM} \z^{-ij} \frac{x^{d-j}-\z^{ij} y^{d-j}}{x-\z^{-i}y} s_i 
\Bigl[x,\frac{X^j-\z^{ij}Y^j}{X-\z^i Y}\Bigr] } \\
&=& \DS{a^2 \sum_{i,i' \in \ZM/d\ZM} 
\z^{-ij} \frac{x^{d-j}-\z^{ij} y^{d-j}}{x-\z^{-i}y} s_i s_{i'} \Bigl(\frac{X^j-\z^{ij}Y^j}{X-\z^i Y}-
\frac{\z^{i'j}X^j-\z^{(i-i')j} Y^j}{\z^{i'}X-\z^{i-i'} Y} \Bigr)},
\endeqna
again by using~(\ref{eq:x-com-p}). 
But $s_is_{i'}=c^{i-i'}$, where $c=ts=\diag(\z,\z^{-1})$ so, if we set $k=i-i'$, 
we can rewrite the above formula as follows:
\eqna
[x,\bbb_j]&=&\DS{a^2 \sum_{i,k \in \ZM/d\ZM} \z^{-ij} \frac{x^{d-j}-\z^{ij} y^{d-j}}{x-\z^{-i}y} c^k 
\Bigl(\frac{X^j-\z^{ij}Y^j}{X-\z^i Y}-
\frac{\z^{(i-k)j}X^j-\z^{kj} Y^j}{\z^{i-k}X-\z^{k} Y} \Bigr)}\\
&=&\DS{a^2 \sum_{i,k \in \ZM/d\ZM} \z^{-ij} \frac{x^{d-j}-\z^{ij} y^{d-j}}{x-\z^{-i}y} 
\Bigl(\frac{\z^{-kj} X^j-\z^{(i+k)j}Y^j}{\z^{-k}X-\z^{i+k} Y}-
\frac{X^j-\z^{ij} Y^j}{X-\z^i Y} \Bigr) c^k}\\
&=&\DS{a^2 \sum_{k \in \ZM/d\ZM} \Th_{j,k} c^k,}
\endeqna
where 
$$\Th_{j,k}=\sum_{i \in \ZM/d\ZM} 
\z^{-ij} \frac{x^{d-j}-\z^{ij} y^{d-j}}{x-\z^{-i}y} 
\Bigl(\frac{\z^{-kj} X^j-\z^{(i+k)j}Y^j}{\z^{-k}X-\z^{i+k} Y}-
\frac{X^j-\z^{ij} Y^j}{X-\z^i Y}\Bigr) \in \CM[x,y] \otimes \CM[X,Y].$$
This formula implies that $\Th_{j,k}$ is a linear combination of (non-commutative) 
monomials of the form $x^ly^{d-1-l}X^m Y^{j-1-m}$, 
where $0 \le l \le d-j-1$ and $0 \le m \le j-1$, and the coefficient $\th_{j,k,l,m}$ 
of this monomial in $\Th_{j,k}$ is equal to 
\eqna
\th_{j,k,l,m}&=&\DS{\sum_{i \in \ZM/d\ZM} \z^{-ij}\z^{-i(d-j-1-l)}(\z^{-km}\z^{(i+k)(j-1-m)} 
- \z^{i(j-1-m)})}\\
&=& \DS{\sum_{i \in \ZM/d\ZM} \z^{i(l+j-m)}(\z^{k(j-1-2m)}-1)}.
\endeqna
But $j \le l+j \le d-1$ and $0 \le m \le j-1$, so $l+j \not\equiv m \mod d$. 
This implies in particular that $\sum_{i \in \ZM/d\ZM} \z^{i(l+j-m)}=0$, and so $\th_{j,k,l,m}=0$. 
This shows that $[x,\bbb_j]=0$.

A similar computation shows that $[X,\bbb_j]=0$ and so $\bbb_j$ commutes with $s$, $t$, 
$x$, $X$, $sxs^{-1}=y$ and $sXs^{-1}=Y$, so it is central in $\Hb_c$. This completes
the proof of the proposition.
\end{proof}

\bigskip

This has the following consequence, that will be used later for obtaining 
a presentation of the algebra $Z_c$.

\medskip

\begin{coro}\label{cor:horreur}
If $1 \le i \le j \le d-1$, then 
\begin{multline*}
\troncation_c(\abb_{i-1}\abb_{j+1}-\abb_i\abb_j) = 
q^{d-j-1} Q^{i-1} (x^{j-i+2} X^{j-i+2}+y^{j-i+2}Y^{j-i+2}) \\
-q^{d-j} Q^i (x^{j-i} X^{j-i}+y^{j-i}Y^{j-i}) \\
\DS{+ d(1+j-i-d) a^2 \sum_{M=i-1}^{j-1}    
x^{M+d-i-j}y^{d-2-M} X^M Y^{i+j-2-M}}. 
\end{multline*}
\end{coro}

\medskip

\begin{proof}
Assume first that $0 \le i \le j \le d$. 
Since 
$\abb_j$ is central, we get 
\eqna
\abb_i\abb_j &=& x^{d-i} \abb_j Y^i + y^{d-i} \abb_j X^i 
\DS{- a \sum_{k \in \ZM/d\ZM} \z^{-ki} \g_{k,i} \abb_j s_k \G_{k,i} } \\
&=& x^{d-i}x^{d-j}Y^j Y^i + x^{d-i} y^{d-j} X^j Y^i 
+ y^{d-i} x^{d-j} Y^j X^i + y^{d-i} y^{d-j} X^j X^i \\
&& \DS{-a \sum_{k \in \ZM/d\ZM} x^{d-i} \g_{k,j} s_k \G_{k,j} Y^i 
-a \sum_{k \in \ZM/d\ZM} y^{d-i} \g_{k,j} s_k \G_{k,j} X^i} \\
&&\DS{-a \sum_{k \in \ZM/d\ZM} \z^{-ki} \g_{k,i} 
(x^{d-j}Y^j+y^{d-j}X^j) \G_{k,i} s_k} \\
&&\DS{+a^2 \sum_{k,l \in \ZM/d\ZM} \z^{-ki}\z^{-kj} 
\g_{k,i} \g_{l,j} \G_{l,j} s_l s_k \G_{k,i}} \\
\endeqna
Therefore, 
\eqna
\troncation_c(\abb_i\abb_j) &=& x^{2d-i-j}Y^{i+j} + y^{2d-i-j}X^{i+j} 
+ q^{d-j} Q^i (x^{j-i} X^{j-i}+y^{j-i}Y^{j-i}) \\
&& \DS{+ a^2 \sum_{k \in \ZM/d\ZM} \z^{-k(i+j)}
\g_{k,i} \g_{k,j} \G_{k,j} \G_{k,i}}.
\endeqna
Expanding the product $\g_{k,i} \g_{k,j} \G_{k,j} \G_{k,i}$ gives
\begin{multline*}
\troncation_c(\abb_i\abb_j) = x^{2d-i-j}Y^{i+j} + y^{2d-i-j}X^{i+j} 
+ q^{d-j} Q^i (x^{j-i} X^{j-i}+y^{j-i}Y^{j-i}) \\
\DS{+ a^2 \sum_{k \in \ZM/d\ZM} \sum_{l=0}^{d-i-1}~\sum_{l'=0}^{d-j-1}~ 
\sum_{m=0}^{i-1}~ \sum_{m'=0}^{j-1} 
\z^{-k(i+j+l+l'-m-m')} x^{l+l'}y^{2d-i-j-2-l-l'} X^{m+m'} Y^{i+j-2-m-m'}}.
\end{multline*}
If $0 \le L \le 2d-i-j-2$ (resp. $0 \le M \le i+j-2$), let 
$\LC_{i,j}(L)$ (resp. $\MC_{i,j}(M)$) denote the set of pairs 
$(l,l')$ (resp. $(m,m')$) such that $l+l'=L$ (resp. $m+m'=M$) 
and $0 \le l \le d-i-1$ and $0 \le l' \le d-j-1$ (resp. 
$0 \le m \le i-1$ and $0 \le m' \le j-1$). Then the above equality 
might rewritten
\begin{multline*}
\troncation_c(\abb_i\abb_j) = x^{2d-i-j}Y^{i+j} + y^{2d-i-j}X^{i+j} 
+ q^{d-j} Q^i (x^{j-i} X^{j-i}+y^{j-i}Y^{j-i}) \\
\DS{+ a^2 \sum_{k \in \ZM/d\ZM} \sum_{L=0}^{2d-i-j-2} \hphantom{A}
\sum_{M=0}^{i+j-2}  
|\LC_{i,j}(L)|\cdot |\MC_{i,j}(M)| \cdot \z^{-k(i+j+L-M)} 
x^Ly^{2d-i-j-2-L} X^M Y^{i+j-2-M}}. 
\end{multline*}
Now, if $1 \le i \le j \le d-1$, applying the above formula 
by replacing $i$ by $i-1$ and $j$ by $j+1$ yields 
\begin{multline*}
\troncation_c(\abb_{i-1}\abb_{j+1}-\abb_i\abb_j) = 
q^{d-j-1} Q^{i-1} (x^{j-i+2} X^{j-i+2}+y^{j-i+2}Y^{j-i+2}) \\
-q^{d-j} Q^i (x^{j-i} X^{j-i}+y^{j-i}Y^{j-i}) \\
\DS{+ a^2 \sum_{L=0}^{2d-i-j-2} \hphantom{A}
\sum_{M=0}^{i+j-2}  \Bigl(\sum_{k \in \ZM/d\ZM} \z^{-k(i+j+L-M)} \Bigr) 
(|\LC_{i-1,j+1}(L)|\cdot |\MC_{i-1,j+1}(M)|-|\LC_{i,j}(L)|\cdot |\MC_{i,j}(M)|)} \\ 
x^Ly^{2d-i-j-2-L} X^M Y^{i+j-2-M}. 
\end{multline*}
So, the coefficient of $x^Ly^{2d-i-j-2-L} X^M Y^{i+j-2-M}$ 
is non-zero if and only if $i+j+L \equiv M \mod d$ and 
$|\LC_{i-1,j+1}(L)|\cdot |\MC_{i-1,j+1}(M)|\neq|\LC_{i,j}(L)|\cdot |\MC_{i,j}(M)|$. 
Since $i \le j$, we have
$$|\LC_{i,j}(L)|=
\begin{cases}
1+L & \text{if $0 \le L \le d-j-1$,}\\
d-j & \text{if $d-j-1 \le L \le d-i-1$,}\\
2d-i-j-1-L & \text{if $d-i-1 \le L \le 2d-i-j-2$,}\\
\end{cases}$$
$$|\MC_{i,j}(M)|=
\begin{cases}
1+M & \text{if $0 \le M \le i-1$,}\\
i & \text{if $i-1 \le M \le j-1$,}\\
i+j-1-M & \text{if $j-1 \le M \le i+j-2$.}
\end{cases}\leqno{\text{and}}
$$
So $|\LC_{i-1,j+1}(L)|\cdot |\MC_{i-1,j+1}(M)|\neq|\LC_{i,j}(L)|\cdot |\MC_{i,j}(M)|$ 
if and only if $d-j-1 \le L \le d-i-1$ or $i-1 \le M \le j-1$. Combined 
with the fact that $i+j+L \equiv M \mod d$ to obtain a non-zero 
coefficient for $x^Ly^{2d-i-j-2-L} X^M Y^{i+j-2-M}$, this forces $i+j+L=M+d$ 
and so
\begin{multline*}
\troncation_c(\abb_{i-1}\abb_{j+1}-\abb_i\abb_j) = 
q^{d-j-1} Q^{i-1} (x^{j-i+2} X^{j-i+2}+y^{j-i+2}Y^{j-i+2}) \\
-q^{d-j} Q^i (x^{j-i} X^{j-i}+y^{j-i}Y^{j-i}) \\
\DS{+ da^2 \sum_{M=i-1}^{j-1}  
\underbrace{((d-j-1)(i-1)-(d-j)i)}_{=1+j-i-d} x^{M+d-i-j}y^{d-2-M} X^M Y^{i+j-2-M}}, 
\end{multline*}
as expected.
\end{proof}

\subsection{Poisson bracket}
We determine here part of the Poisson bracket between the generators:

\smallskip

\begin{prop}\label{prop:poisson}
We have 
$$\{q,Q\}=\eulerb,\quad \{\eulerb,q\}=-2q\quad\text{\it and}\quad \{\eulerb,Q\}=2Q.$$
Moreover, if $0 \le j \le d$, then
$$\{q,\abb_j\}=j \abb_{j-1},\quad \{\eulerb,\abb_j\}= (2j-d) \abb_j\quad\text{\it and}\quad 
\{Q,\abb_j\}=(j-d) \abb_{j+1},$$
with the convention that $\abb_{-1}=\abb_{d+1}=0$. 
\end{prop}

\medskip

\begin{proof}
First, note that the Poisson bracket on $Z_c$ is in 
fact the restriction of a Poisson bracket $\{,\} : \Hb_c \times Z_c \longto \Hb_c$. 
This Poisson bracket satisfies the following property: if $z=\sum_{w \in W} f_w w F_w \in Z_c$, 
with $f_w \in \CM[x,y]$ and $F_w \in \CM[X,Y]$, then 
\equat\label{eq:poisson-x}
\{x,z\}=\sum_{w \in W} f_w w \frac{\partial F_w}{\partial X}
\qquad\text{and}\qquad
\{y,z\}=\sum_{w \in W} f_w w \frac{\partial F_w}{\partial Y}.
\endequat
The first three equalities of the proposition are standard and hold for any Coxeter group 
(see~\cite[\S{4}]{dezelee} or~\cite[\S{3}]{berest}) and can easily be checked 
in this case by a little computation. 
Similarly, the fact that $\{\eulerb,\abb_j\}= (2j-d) \abb_j$ follows from the general 
fact that, if $h \in \Hb_c$ is homogeneous of degree $k$, then $\{\eulerb,h\}=kh$ 
(see for instance~\cite[Prop.~3.3.3]{calogero}). We now prove that $\{q,\abb_j\}=j \abb_{j+1}$, 
the last equality 
being proved similarly. 
From the formula given for $\abb_j$ in Proposition~\ref{prop:ai}, 
we get
\eqna
\{q,\abb_j\}=\{yx,\abb_j\} &=& j x^{d-j} Y^{j-1} x + j y^{d-j+1} X^{j-1} \\
&&\DS{-a \sum_{i \in \ZM/d\ZM} \z^{-ij} \g_{i,j} s_i \frac{\partial \G_{i,j}}{\partial Y} x 
-a \sum_{i \in \ZM/d\ZM} \z^{-ij} y\g_{i,j} s_i \frac{\partial \G_{i,j}}{\partial X}.}
\endeqna
In order to prove the proposition, it is sufficient to check that 
$\troncation_c(\{q,\abb_j\})=j \ab_{j-1}$. But, from the above formula and from~(\ref{eq:x-com-p}), 
one gets 
$$
\troncation_c(\{q,\abb_j\})= j \ab_{j-1} -a 
\troncation_c\Bigl(\sum_{i \in \ZM/d\ZM} \z^{-ij} \g_{i,j} s_i 
\frac{\partial \G_{i,j}}{\partial Y} x \Bigr).$$
Since
$$\sum_{i \in \ZM/d\ZM} \z^{-ij} \g_{i,j} s_i 
\frac{\partial \G_{i,j}}{\partial Y} x = 
\sum_{i \in \ZM/d\ZM} \z^{-ij} \g_{i,j} s_i 
(x\frac{\partial \G_{i,j}}{\partial Y} - \Bigl[x,\frac{\partial \G_{i,j}}{\partial Y}\Bigr]),$$
it follows from~(\ref{eq:x-com-p}) that
$$\troncation_c(\{q,\abb_j\})= j \ab_{j-1} + 
a^2 \sum_{i \in \ZM/d\ZM} \z^{-ij} \g_{i,j} \frac{\frac{\partial \G_{i,j}}{\partial Y} 
- \lexp{s_i}{\Bigl(\frac{\partial \G_{i,j}}{\partial Y}\Bigr)}}{X-\z^i Y}.$$
So it remains to prove that 
$$\sum_{i \in \ZM/d\ZM} \z^{-ij} \g_{i,j} \frac{\frac{\partial \G_{i,j}}{\partial Y} 
- \lexp{s_i}{\Bigl(\frac{\partial \G_{i,j}}{\partial Y}\Bigr)}}{X-\z^i Y}=0.\leqno{(?)}$$
Let us compute the big fraction in the above formula. First,
$$\G_{i,j}=\sum_{k=0}^{j-1} \z^{ik} X^{j-1-k} Y^k,$$
so
$$\frac{\partial \G_{i,j}}{\partial Y}=\sum_{k=0}^{j-1} k\z^{ik} X^{j-1-k} Y^{k-1}
=\sum_{k=0}^{j-2} (k+1)\z^{i(k+1)} X^{j-2-k} Y^{k}.$$
Therefore,
$$\lexp{s_i}{\Bigl(\frac{\partial \G_{i,j}}{\partial Y}\Bigr)}
=\sum_{k=0}^{j-2} (k+1)\z^{i(k+1)}(\z^i Y)^{j-2-k} (\z^{-i}X)^{k}.$$
Simplifying and using the change of variable $k \mapsto j-2-k$, one gets
$$\lexp{s_i}{\Bigl(\frac{\partial \G_{i,j}}{\partial Y}\Bigr)}
=\sum_{k=0}^{j-2} (j-1-k)\z^{i(k+1)}X^{j-2-k} Y^{k}.$$
We deduce that
$$\frac{\partial \G_{i,j}}{\partial Y}-\lexp{s_i}{\Bigl(\frac{\partial \G_{i,j}}{\partial Y}\Bigr)}
= \sum_{k=0}^{j-2} (2k+2-j) \z^{i(k+1)} X^{j-2-k} Y^{k}.$$
But $\sum_{k=0}^{j-2} (2k+2-j) = 0$, so 
$$\frac{\partial \G_{i,j}}{\partial Y}-\lexp{s_i}{\Bigl(\frac{\partial \G_{i,j}}{\partial Y}\Bigr)}
= \sum_{k=0}^{j-2} (2k+2-j) \z^{i(k+1)} (X^{j-2-k} - \z^{i(j-2-k)}Y^{j-2-k})Y^{k}.$$
Since the term corresponding to $k=j-2$ vanishes, this implies that
\eqna
\frac{\frac{\partial \G_{i,j}}{\partial Y} 
- \lexp{s_i}{\Bigl(\frac{\partial \G_{i,j}}{\partial Y}\Bigr)}}{X-\z^i Y}
&=& \DS{\sum_{k=0}^{j-3} \sum_{k'=0}^{j-3-k} (2k+2-j)\z^{i(k+1)} X^{j-3-k-k'}(\z^i Y)^{k'} Y^k}\\
&=& \DS{\sum_{k=0}^{j-3} \sum_{k'=0}^{j-3-k} (2k+2-j)\z^{i(k+k'+1)} X^{j-3-k-k'} Y^{k+k'}}\\
&=& \DS{\sum_{k=0}^{j-3} \sum_{k'=k}^{j-3} (2k+2-j)\z^{i(k'+1)} X^{j-3-k'} Y^{k'}}\\
&=& \DS{\sum_{k'=0}^{j-3} \Bigl(\sum_{k=0}^{k'} (2k+2-j)\Bigr)\z^{i(k'+1)} X^{j-3-k'} Y^{k'}}\\
&=& \DS{\sum_{k'=0}^{j-3} (k'+2-j)(k'+1) \z^{i(k'+1)} X^{j-3-k'} Y^{k'}}.
\endeqna
Therefore, the left-hand side of the formula~(?) is a linear combination 
of monomials of the form $x^{d-j-1-l}y^lX^{j-3-m}Y^m$, where $0 \le l \le d-j-1$ 
and $0 \le m \le j-3$, and the coefficient of this monomial is
$$\sum_{i \in \ZM/d\ZM} (m+2-j)(m+1) \z^{-ij}\z^{-il}\z^{i(m+1)}=(m+2-j)(m+1)\sum_{i \in \ZM/d\ZM} 
\z^{i(m+1-j-l)}.$$
But $j \le l+j \le d-1$ and $1 \le m+1 \le j-2$, 
so this coefficient is $0$ and the equality~(?) is proved.
\end{proof}

\bigskip

\begin{coro}\label{cor:invariant}
We have
$$\{q,\eulerb^2-4qQ\}=\{\eulerb,\eulerb^2-4qQ\}=\{Q,\eulerb^2-4qQ\}=0.$$
\end{coro}

\bigskip

\subsection{Presentation} 
The main result of this paper is the following:

\bigskip

\begin{theo}\label{theo:presentation-zc}
If $a=b$, then the algebra $Z_c$ admits the following presentation:
\begin{itemize}
\item[$\bullet$] Generators: $q$, $Q$, $\eulerb$, $\abb_0$, $\abb_1$, \dots, $\abb_d$.

\item[$\bullet$] Relations:
$$
\begin{cases}
\eulerb \,\abb_i = q \abb_{i+1} + Q \abb_{i-1} & \text{\it for $1 \le i \le d-1$,}\\
\abb_{i-1} \abb_{j+1} - \abb_i \abb_j  = 
(\eulerb^2-4qQ-d^2 a^2) q^{d-j-1} Q^{i-1} \Psi_{j-i}(\eulerb,q,Q)  
& \text{\it for $1 \le i \le j \le d-1$.}\\
\end{cases}
$$
\end{itemize}
\end{theo}

\bigskip

\begin{proof}
By~\cite{bonnafe thiel}, a presentation of $Z_c$ is obtained by deforming 
the generators of $Z_0=\CM[V \times V^*]^W$ and deforming the relations. Therefore, 
in order to prove the theorem, it is sufficient to check that the relations 
given in the statement are satisfied. So let $1 \le i \le j \le d-1$.

\medskip

Let us first prove that 
$$\eulerb \,\abb_i = q \abb_{i+1} + Q \abb_{i-1}.\leqno{(\ZG_i)}$$
For this, it is sufficient to prove that 
$\troncation_c(\eulerb \,\abb_i) = \troncation_c(q \abb_{i+1} + Q \abb_{i-1})$. 
But the map $\troncation_c$ is $P_{\!\bullet}$-linear 
so it is sufficient to prove that
$$\troncation_c(\eulerb \,\abb_i)=\euler_0\, \ab_{i,0}.\leqno{(?)}$$
Since $\abb_i$ is central, it follows from~(\ref{eq:euler}) 
and Proposition~\ref{prop:ai} that
\eqna
\eulerb \,\abb_i &=& x \abb_i X + y \abb_i Y + a \DS{\sum_{k \in \ZM/d\ZM} \abb_i s_k}\\
&=& x (x^{d-i}Y^i+y^{d-i}X^i)X + y (x^{d-i}Y^i+y^{d-i}X^i) Y \\
&& \DS{- a \sum_{k' \in \ZM/d\ZM} \z^{-ik'} x \g_{k',i} s_{k'} \G_{k',i} X - 
a \sum_{k' \in \ZM/d\ZM} \z^{-ik'} y \g_{k',i} s_{k'} \G_{k',i} Y} \\
&& + a \sum_{k \in \ZM/d\ZM} (x^{d-i}Y^i+y^{d-i}X^i) s_k 
- a^2 \sum_{k,k' \in \ZM/d\ZM} \z^{-ik'} \g_{k',i}\G_{k',i} s_{k'} s_k.
\endeqna
But $s_{k'}s_k=1$ if and only if $k'=k$, so 
$$\troncation_c(\eulerb \,\abb_i)=\euler_0 \ab_{i,0} - a^2 
\sum_{k \in \ZM/d\ZM} \z^{-ik} \g_{k,i}\G_{k,i}.$$
The element $\sum_{k \in \ZM/d\ZM} \z^{-ik} \g_{k,i}\G_{k,i}$ of $\CM[V \times V^*]^W$ is a linear 
combination of monomials of the form $x^{d-i-1-l}y^l X^{i-1-m}Y^m$ 
where $0 \le l \le d-i-1$ and $0 \le m \le i-1$, and the coefficient 
of this monomial is equal to
$$\sum_{k \in \ZM/d\ZM} \z^{-ki} \z^{-kl} \z^{km} = \sum_{k \in \ZM/d\ZM} 
\z^{-k(i+l-m)}.$$
But $i \le i+l \le d-1$ and $0 \le m \le i-1$, so $i+l-m \not\equiv 0 \mod d$. 
This shows that the above sum is zero, and this completes the proof of~(?).

\medskip

Let us now prove that 
$$\abb_{i-1}\abb_{j+1} - \abb_i\abb_j=
(\eulerb^2-4qQ-d^2 a^2) q^{d-j-1} Q^{i-1} \Psi_{j-i}(\eulerb,q,Q).\leqno{(\ZG_{i,j})}$$
This will be proved by induction on $j-i$. So let us first consider the case 
where $j-i=0$, i.e. where $j=i$. 
Again, it is sufficient to prove the equality after applying the map 
$\troncation_c$. We deduce from 
Corollary~\ref{cor:horreur} that 
$$\troncation_c(\abb_{i-1}\abb_{i+1}-\abb_i^2) = 
q^{d-i-1} Q^{i-1} (x^2 X^2+y^2Y^2-2qQ - d(d-1)a^2).$$
Since $\troncation_c$ is $P_{\!\bullet}$-linear and $\Psi_0=1$, 
proving $(\ZG_{i,i})$ is equivalent to proving that 
$$\troncation_c(\eulerb^2-4qQ-d^2 a^2)=x^2 X^2+y^2Y^2-2qQ - d(d-1)a^2,$$
or, equivalently, that 
$$\troncation_c(\eulerb^2)=x^2 X^2+y^2Y^2+2qQ +d a^2.\leqno{(\EC)}$$
But 
\eqna
\eulerb^2&=& x\eulerb X + y \eulerb Y + a \DS{\sum_{k \in \ZM/d\ZM} \eulerb s_k} \\
&=&\DS{x^2X^2+xyYX + a \sum_{k \in \ZM/d\ZM} xs_kX +yxXY+ y^2Y^2 + 
a \sum_{k \in \ZM/d\ZM} ys_kY }\\
&& \DS{+ a \sum_{k \i \ZM/d\ZM} (xX+yY) s_k + a^2 \sum_{k,l \in \ZM/d\ZM} s_ls_k.}
\endeqna
It follows directly that 
$$\troncation_c(\eulerb^2)=x^2 X^2+y^2Y^2+2qQ + a^2 \sum_{k \in \ZM/d\ZM} 1,$$
as desired.

\medskip

Assume now that $j-i \ge 1$ and that $(\ZG_{i',j'})$ holds if $j'-i' < j-i$. 
Then, by the induction hypothesis, we have
$$\abb_i\abb_{j+1}-\abb_{i+1}\abb_j=
(\eulerb^2-4qQ-d^2 a^2) q^{d-j-2} Q^i\Psi_{j-i-1}(\eulerb,q,Q).$$
Applying $\{q,-\}$ to this equality, and using Proposition~\ref{prop:poisson} 
and Corollary~\ref{cor:invariant}, one gets:
\begin{multline*}
i(\abb_{i-1}\abb_{j+1}-\abb_i\abb_j) + j (\abb_i\abb_j-\abb_{i+1}\abb_{j-1})= 
(\eulerb^2-4qQ-d^2 a^2) q^{d-j-1} \\
\times \bigl(iQ^{i-1}\eulerb \Psi_{j-i-1}(\eulerb,q,Q) 
+ 2q Q^i \DS{\frac{\partial \Psi_{j-i-1}}{\partial T}(\eulerb,q,Q)
+ Q^i \eulerb \frac{\partial \Psi_{j-i-1}}{\partial T''}(\eulerb,q,Q) \bigr)}.
\end{multline*}
But, by~(\ref{eq:pi-der}), 
$$2q \frac{\partial \Psi_{j-i-1}}{\partial T}(\eulerb,q,Q) + 
\eulerb \frac{\partial \Psi_{j-i-1}}{\partial T''}(\eulerb,q,Q) =(j-i) q\Psi_{j-i-2}(\eulerb,q,Q)$$
and, by~(\ref{eq:pi-rec}), 
$$\eulerb \Psi_{j-i-1}(\eulerb,q,Q) - qQ \Psi_{j-i-2}(\eulerb,q,Q)=\Psi_{j-i}(\eulerb,q,Q).$$
Therefore,
\begin{multline*}
i(\abb_{i-1}\abb_{j+1}-\abb_i\abb_j) + j (\abb_i\abb_j-\abb_{i+1}\abb_{j-1})=  \\
(\eulerb^2-4qQ-d^2 a^2) q^{d-j-1}Q^{j-1}(i\Psi_{j-i}(\eulerb,q,Q) + j qQ \Psi_{j-i-2}(\eulerb,q,Q)).
\end{multline*}
Since the induction hypothesis implies that 
$$(\abb_i\abb_j-\abb_{i+1}\abb_{j-1})= (\eulerb^2-4qQ-d^2 a^2) q^{d-j}Q^i\Psi_{j-i-2}(\eulerb,q,Q),$$
the result follows.
\end{proof}

\bigskip

\subsection{Back to Poisson bracket}
In Proposition~\ref{prop:poisson}, we did not determine the Poisson brackets 
$\{\abb_i,\abb_j\}$. This was only determined for $a=0$ in~(\ref{eq:poisson-0}): 
it is proven that there exists a polynomial $\Pi_{i,j} \in \CM[T,T',T'']$, 
which is homogeneous of degree $d-1$, such that
$$\{\ab_{i,0},\ab_{j,0}\}=\Pi_{i,j}(\euler_0,q,Q).$$
This will be deformed to the unequal parameter case as follows:

\bigskip

\begin{prop}\label{prop:poisson-ai-aj}
If $0 \le i < j \le d$, there exists a polynomial 
$\Phi_{i,j} \in \CM[T,T',T'']$, homogeneous of degree $d-3$, such that
$$\{\abb_i,\abb_j\}=\Pi_{i,j}(\eulerb,q,Q)+a^2\Phi_{i,j}(\eulerb,q,Q).$$
\end{prop}

\bigskip

\begin{proof}
We will prove that there exist polynomials 
$\Phi_{i,j}^\circ$, $\Phi_{i,j} \in \CM[T,T',T'']$, homogeneous of degree $d-1$ and 
$d-3$ respectively, such that 
$$\{\abb_i,\abb_j\}=\Phi_{i,j}^\circ(\eulerb,q,Q)+a^2\Phi_{i,j}(\eulerb,q,Q).\leqno{(\wp_{i,j})}$$
This is sufficient because, by specializing $a$ to $0$, one gets that $\Phi_{i,j}^\circ=\Pi_{i,j}$. 

Let us first assume that $i=0$. To make an induction argument on $j$ work, 
we will prove a slightly stronger result, namely that 
$$
\{\abb_0,\abb_j\}=q^{d-j}(\ph_j(\eulerb,q,Q)+a^2\th_j(\eulerb,q,Q)).\leqno{(\wp_{0,j}^+)}
$$
where $\ph_j$, $\th_j \in \CM[T,T',T'']$ are homogeneous of degree $j-1$ and 
$j-3$ respectively.
For this, let us apply $\{\abb_0,-\}$ to the following two relations given by 
Theorem~\ref{theo:presentation-zc}
$$
Q \abb_0 -\eulerb \abb_1 + q \abb_2=0,\leqno{(\ZG_1)}
$$
$$
\abb_0 \abb_2-\abb_1^2 = q^{d-1}(\eulerb^2 - 4 qQ - d^2a^2).\leqno{(\ZG_{1,1})}
$$
Using Proposition~\ref{prop:poisson}, this gives
$$
\begin{cases}
-\eulerb \{\abb_0,\abb_1\} + q \{\abb_0,\abb_2\} =0,\\
q\abb_0 \{\abb_0,\abb_2\}-2q\abb_1 \{\abb_0,\abb_1\} 
= q^{d-1} \{\abb_0,\eulerb^2-4qQ\}=q^{d-1}(2d\abb_0\eulerb-4dq\abb_1).
\end{cases}
$$
Thanks to the first equality, we can replace the term $q \{\abb_0,\abb_2\}$ in the 
second equation by $\eulerb \{\abb_0,\abb_1\}$, and this yields
$$(\abb_0\eulerb-2q\abb_1) \{\abb_0,\abb_1\}= 2dq^{d-1}(\abb_0\eulerb-2q\abb_1).$$
Since $\abb_0\eulerb-2q\abb_1 \neq 0$ (by computing its image by $\troncation_c$) 
and since $Z_c$ is an integral domain, we get 
$$\{\abb_0,\abb_1\}= 2dq^{d-1},$$
which proves $(\wp_{0,1}^+)$. We also deduce that 
$$\{\abb_0,\abb_2\}= 2dq^{d-2}\eulerb,$$
which proves $(\wp_{0,2}^+)$. 

Now, assume that $j \ge 3$ and that $(\wp_{0,j'}^+)$ holds for $j' < j$. Applying 
$\{\abb_0,-\}$ to 
$$Q \abb_{j-2}-\eulerb \abb_{j-1} + q \abb_j=0\leqno{(\ZG_{j-1})}$$
yields, thanks to Proposition~\ref{prop:poisson}, 
$$d \abb_1\abb_{j-2}+Q \{\abb_0,\abb_{j-2}\} - d \abb_0 \abb_{j-1} - \eulerb \{\abb_0,\abb_{j-1}\}
+ q \{\abb_0,\abb_j\}=0.$$
But
$$\abb_0\abb_{j-1}-\abb_1\abb_{j-2} = (\eulerb^2-4qQ -d^2a^2) q^{d-j+1} \Psi_{j-3}(\eulerb,q,Q)
\leqno{(\ZG_{1,j-2})}$$
by Theorem~\ref{theo:presentation-zc} and
$$\{\abb_0,\abb_{j-2}\}=q^{d-j+2}(\ph_{j-2}(\eulerb,q,Q)+a^2\th_{j-2}(\eulerb,q,Q)),$$
$$\{\abb_0,\abb_{j-1}\}=q^{d-j+1}(\ph_{j-1}(\eulerb,q,Q)+a^2\th_{j-1}(\eulerb,q,Q))$$
by the induction hypothesis. This gives
\begin{multline*}
\{\abb_0,\abb_j\}=d(\eulerb^2-4qQ -d^2a^2) q^{d-j} \Psi_{j-3}(\eulerb,q,Q) \\
- q^{d-j+1} Q (\ph_{j-2}(\eulerb,q,Q)+a^2\th_{j-2}(\eulerb,q,Q)) + q^{d-j}\eulerb 
(\ph_{j-1}(\eulerb,q,Q)+a^2\th_{j-1}(\eulerb,q,Q)),
\end{multline*}
which proves that $(\wp_{0,j}^+)$ holds.

\medskip

We will now prove that $(\wp_{i,j})$ holds by induction on $i$. The case $i=0$ 
has just been treated, so assume that $i \ge 1$ and that $(\wp_{i-1,j'})$ holds 
for all $j'$. Then $(i-1-d) \abb_i= \{Q,\abb_{i-1}\}$ and $i-1-d \neq 0$. 
By the Jacobi identity, we get
\eqna
(i-1-d) \{\abb_i,\abb_j\}&=&\{\{Q,\abb_{i-1}\},\abb_j\} \\
&=& \{Q,\{\abb_{i-1},\abb_j\}\} - \{\abb_{i-1},\{Q,\abb_j\}\}\\
&=& \{Q,\{\abb_{i-1},\abb_j\}\} - (j-d) \{\abb_{i-1},\abb_{j+1}\}.
\endeqna
So the result follows from the induction hypothesis because, if $\Th \in \CM[T,T',T'']$ 
is an homogeneous polynomial of degree $k$, then 
$$\{Q,\Th(\eulerb,q,Q)\}=-2Q \frac{\partial \Th}{\partial T}(\eulerb,q,Q) - 
\eulerb \frac{\partial \Th}{\partial T'}(\eulerb,q,Q)$$
is of the form $\Th^\#(\eulerb,q,Q)$ where $\Th^\#$ is homogeneous of degree $k$. 
The proof of the proposition is complete.
\end{proof}

\bigskip

\subsection{Lie algebra structure at the cuspidal point}\label{sub:lie}
By Theorem~\ref{theo:presentation-zc}, the affine variety $\ZC_c$ might be described as
\begin{multline*}
\ZC_c=\{(\qG,\QG,\eG,a_0,a_1,\dots,a_d) \in \CM^{d+4}~|~\\
\forall~1\le i \le j \le d-1,~
\begin{cases}
\eG a_i = \qG a_{i+1} + \QG a_{i-1},\\
a_{i-1} a_{j+1} - a_i a_j  = 
(\eG^2-4\qG\QG-d^2 a^2) \qG^{d-j-1} \QG^{i-1} \Psi_{j-i}(\eG,\qG,\QG)
\end{cases} \}.
\end{multline*}
If $d=3$ and $a \neq 0$, then $\ZC_c$ is smooth. So assume from now on that $d \ge 4$ and $a \neq 0$. 
Then the homogeneous component of minimal degree of all 
the above equations is equal to $2$, so the point $0=(0,...,0) \in \ZC_c$ is singular 
and the tangent space of $\ZC_c$ at $0$ has dimension $d+4$. 
It is the only singular point and it is a cuspidal 
point in the sense of~\cite{bellamy cuspidal} (see~\cite[\S{5.2}]{bonnafe diedral}). 
This means that the corresponding maximal ideal $\mG_0$ of $Z_c$ 
is a Poisson ideal (since $\mG_0=\langle q,Q,\eulerb,\abb_0,\abb_1,\dots,\abb_d\rangle$, 
this can also be checked thanks to Proposition~\ref{prop:poisson-ai-aj}). 
This implies that the cotangent space $\mG_0/\mG_0^2$ of $\ZC_c$ at $0$ inherits a Lie algebra structure 
from the Poisson bracket: we denote by $\Lie_0(\ZC_c)$ 
the vector space $\mG_0/\mG_0^2$ endowed with its Lie algebra structure. 
It has been proved in~\cite[Prop.~8.4]{bonnafe diedral} that
\equat\label{eq:lie}
\text{\it If $d=4$, then $\Lie_0(\ZC_c) \simeq \sG\lG_3(\CM)$.}
\endequat
We now determine $\Lie_0(\ZC_c)$ in the remaining cases:

\bigskip

\begin{prop}\label{prop:lie}
If $d \ge 5$, then 
$$\Lie_0(\ZC_c) = \sG\lG_2(\CM) \oplus S_d,$$
where $S_d$ is a commutative ideal of $\Lie_0(\ZC_c)$ of dimension $d+1$ on 
which $\sG\lG_2(\CM)$ acts irreducibly (i.e. $S_d \simeq \Srm\yrm\mrm^d(\CM^2)$ 
as an $\sG\lG_2(\CM)$-module).
\end{prop}

\bigskip

\begin{proof}
If $m \in \mG_0$, we denote by $\mdo$ its image in $\Lie_0(\ZC_c)$. 
Then $(\dot{q},\dot{Q},\dot{\eulerb},\dot{\abb}_0,\dot{\abb}_1,\dots,\dot{\abb}_d)$ 
is a basis of $\Lie_0(\ZC_c)$. We set 
$$\gG=\CM \dot{Q} \oplus \CM \dot{\eulerb} \oplus \CM \dot{q}
\quad\text{and}\quad 
S_d = \bigoplus_{j=0}^d \CM \dot{\abb}_j.
$$
It follows from Proposition~\ref{prop:poisson} that $\gG$ is a Lie subalgebra 
of $\Lie_0(\ZC_c)$ isomorphic to $\sG\lG_2(\CM)$, and that $S_d$ is normalized by $\gG$ 
and is isomorphic to $\Srm\yrm\mrm^d(\CM^2)$ as an $\sG\lG_2(\CM)$-module. 

Since $d \ge 5$ (and so $d-3 \ge 2$), 
we get from Proposition~\ref{prop:poisson-ai-aj} that $\{\abb_i,\abb_j\} \in \mG_0^2$ 
and so $[\dot{\abb}_i,\dot{\abb}_j]=0$. This completes the proof of the 
proposition.
\end{proof}

\section{Action of $\Sb\Lb_2(\CM)$}

\medskip

\subsection{Action and Poisson structure} 
The action of $\Sb\Lb_2(\CM)$ on $\CM[V \times V^*] \rtimes W$ deforms 
to an action on $\Hb_c$ by automorphisms of algebras 
as explained for instance in~\cite[\S{3.6}]{calogero}. 
This action commutes with $W$ and is given on elements of $V$ and $V^*$ 
by the same formula as in~\S\ref{sub:sl2-0}. This induces an action of the Lie 
algebra $\sG\lG_2(\CM)$ on $\Hb_c$ by derivations: as in~\S\ref{sub:sl2-0}, 
if $\xi \in \sG\lG_2(\CM)$ and $\ph \in \Hb_c$, we denote by $\xi \bull h$ 
the action of $-\lexp{t}{\xi}$ on $h$. It is related to the Poisson bracket 
through the same formulas as in~\S\ref{sub:sl2-0}:
\equat\label{eq:action-sl2-c}
e \bull \ph = \{Q,\ph\},\qquad h \bull \ph = \{\euler_0,\ph\}
\qquad\text{and}\qquad f \bull \ph =\{-q,\ph\}.
\endequat

\medskip
\def\Sym{{\mathrm{Sym}}}

\subsection{Map to ${\boldsymbol{\sG\lG_2(\CM)}}$}
If $(\qG,\QG,\eG) \in \CM^3$, we denote by $M(\qG,\QG,\eG)$ the matrix
$$M(\qG,\QG,\eG)=\begin{pmatrix} \eG & \QG \\ - \qG & -\eG \end{pmatrix} \in \sG\lG_2(\CM).$$
We identify $\sG\lG_2(\CM)$ with the subspace of $Z_c$ equal to 
$\CM q \oplus \CM Q \oplus \CM \eulerb$ by sending $(e,h,f)$ to $(Q,\eulerb,-q)$: 
by Proposition~\ref{prop:poisson}, 
this identification carries the Lie bracket on $\sG\lG_2(\CM)$ to the Poisson 
bracket on $\CM q \oplus \CM Q \oplus \CM \eulerb$. This gives an identification 
$\CM[q,Q,\eulerb] \simeq \Sym(\sG\lG_2(\CM))$ and the inclusion 
$\CM[q,Q,\eulerb] \subset Z_c$ gives an $\Sb\Lb_2(\CM)$-equivariant 
Poisson map 
$$\mu^* : \ZC_c \longto \sG\lG_2(\CM)^*$$
(the equivariance follows from~(\ref{eq:action-sl2-c})). 
Identifying $\sG\lG_2(\CM)$ with its dual thanks to the trace map 
endows $\sG\lG_2(\CM)$ with a Poisson structure and gives an $\Sb\Lb_2(\CM)$-equivariant 
Poisson map 
$$\mu : \ZC_c \longto \sG\lG_2(\CM).$$ 
The map $\mu$ can be explicitly described by the following formula 
$$
\mu(\qG,\QG,\eG,a_0,a_1,\dots,a_d)=M(\qG,\QG,\eG).
$$

\bigskip

\subsection{Hermite's reciprocity law} 
Let $E = E^\sharp \oplus E_d$ denote the vector space
$$E=\underbrace{\CM Q \oplus \CM \eulerb \oplus \CM q}_{E^\sharp} \oplus 
\underbrace{\CM \abb_0 \oplus \CM \abb_1 \oplus \cdots \oplus \abb_d}_{E_d}.$$
Theorem~\ref{theo:presentation-zc} shows that the natural morphism 
of algebras $\s : \Sym(E) \longto Z_c$ is surjective and it describes its kernel. 
For avoiding the confusion between multiplication in $Z_c$ and 
multiplication in $\Sym(E)$, we will denote by $\star$ the multiplication 
in $\Sym(E)$. For instance, $\abb_0 \star \abb_2 - \abb_1^{\star 2}$ 
is an element of $\Sym(E)$ whereas $\abb_0\abb_2-\abb_1^2$ is an element of $Z_c$, 
which is equal to $\s(\abb_0 \star \abb_2 - \abb_1^{\star 2})$. Similarly, if 
$e_1$,\dots, $e_n$ are elements of $E$ and if $\Psi \in \CM[T_1,\dots,T_n]$ 
is a polynomial in $n$ indeterminates, we denote by $\Psi^\star(e_1,\dots,e_n)$ 
the evaluation of $\Psi$ at $(e_1,\dots,e_n)$ {\it inside the algebra $\Sym(E)$} 
whereas $\Psi(e_1,\dots,e_n)$ denotes the evaluation of $\Psi$ 
{\it inside the algebra $Z_c$}: they satisfy the equality 
$\s(\Phi^\star(e_1,\dots,e_n))=\Psi(e_1,\dots,e_n)$.

Proposition~\ref{prop:poisson} and~(\ref{eq:action-sl2-c}) imply that 
$E$ is an $\Sb\Lb_2(\CM)$-stable subspace of $Z_c$, so that $\s$ is 
$\Sb\Lb_2(\CM)$-equivariant. Let us denote by $V_2 \simeq \CM^2$ another 
copy of $\CM^2$ viewed as the standard representation of $\Sb\Lb_2(\CM)$ 
(or $\sG\lG_2(\CM)$), and we denote by $(t,u)$ its canonical basis. 
We then have two isomorphisms of vector spaces 
$$\s^\sharp : \Sym^2(V_2) \longto E^\sharp \qquad\text{and}\qquad
\s_d : \Sym^d(V_2) \longto E_d$$
which are defined by 
$$\s^\sharp(t^2)=2q,\qquad \s^\sharp(tu)=\eulerb,\qquad \s^\sharp(u^2)=2Q$$
$$\s_d(t^{d-i}u^i)=\abb_i\qquad\text{for $0 \le i \le d$.}\leqno{\text{and}}$$
Proposition~\ref{prop:poisson} and~(\ref{eq:action-sl2-c}) imply that 
$\s^\sharp$ and $\s_d$ are $\Sb\Lb_2(\CM)$-equivariant and we will identify 
$E^\sharp$ and $E_d$ with $\Sym^2(V_2)$ and $\Sym^d(V_2)$ through 
these isomorphisms. 

\def\Der{{\mathrm{Der}}}

Let us first interprete the equations $(\ZG_i)_{1 \le i \le d-1}$. 
Note that 
$$\Sym^2(E)= \Sym^2(\Sym^2(V_2)) ~\oplus~ \Sym^2(V_2) \otimes \Sym^d(V_2) 
~\oplus~ \Sym^2(\Sym^d(V_2))$$ 
and that we have a natural morphism 
$$\mu_{2,d} : \Sym^2(V_2) \otimes \Sym^d(V_2) \longto \Sym^{d+2}(V_2)$$
given by multiplication. We denote by $\Der(\Sym(V_2))$ the $\Sym(V_2)$-module 
of derivations $\Sym(V_2) \to \Sym(V_2)$. If $D \in \Der(\Sym(V_2))$, we denote 
by $D^{(2)}$ the map $\Sym^2(V_2) \otimes \Sym^d(V_2) \longto \Sym(V_2)$, 
$\ph \otimes \psi \longmapsto D(\ph)\psi$. 
Then it is easily checked that 
$$\Ker(\mu_{2,d}) \cap \bigcap_{D \in \Der(\Sym(V_2))} \Ker(D^{(2)}) = \bigoplus_{i=1}^{d-1} 
\CM(Q \star \abb_{i-1} - \eulerb \star \abb_i + q \star \abb_{i+1}) \subset \Sym^2(E).$$
So the family of equations $(\ZG_i)_{1 \le i \le d-1}$ can be summarized by
\equat\label{eq:zi-new}
\text{\it $\Ker(\mu_{2,d})  \cap \bigcap_{D \in \Der(\Sym(V_2))} \Ker(D^{(2)})$ is contained in $\Ker(\s)$.}
\endequat
Note that $\Ker(\mu_{2,d})  \cap \bigcap_{D \in \Der(\Sym(V_2))} \Ker(D^{(2)})$ 
is $\Sb\Lb_2(\CM)$-stable, as the construction is canonical. 

The interpretation of the equations $(\ZG_{i,j})_{1 \le i \le j \le d-1}$ 
is somewhat more subtle and is related with Hermite's reciprocity law (see the upcoming 
Remark~\ref{rem:hermite}). 
First, evaluation induces a surjective morphism of $\Sb\Lb_2(\CM)$-modules 
$$\fonction{\e_{m,n}}{\Sym^m(\Sym^n(V_2))}{\Sym^{mn}(V_2)}{v_1\star \cdots \star v_m}{v_1\cdots v_m.}$$
In the special case where $m=2$ and $n=d$, then: 

\bigskip

\begin{lem}\label{lem:base-ai}
The family $(\abb_{i-1}\star\abb_{j+1}-\abb_i\star\abb_j)_{1 \le i \le j \le d-1}$ of elements 
of $\Sym(E)$ is a basis 
of $\Ker(\e_{2,d}) \subset \Sym^2(\Sym^d(V_2)) \simeq \Sym^2(E_d)$.
\end{lem}

\bigskip

In fact, the family $(\abb_{i-1}\star\abb_{j+1}-\abb_i\star\abb_j)_{1 \le i \le j \le d-1}$ 
generates the ideal equal to the kernel of the natural morphism 
$\e_{\bullet,d} : \Sym(\Sym^d(V_2)) \to \Sym(V_2)$. On the other hand, 
it follows from~(\ref{eq:base}) that:

\bigskip

\begin{lem}\label{lem:base-qq}
The family $(q^{\star d-j-1}\star Q^{\star i-1} \star 
\Psi_{j-i}^\star(\eulerb,q,Q))_{1 \le i \le j \le d-1}$ of elements of $\Sym(E)$ 
is a basis 
of $\Sym^{d-2}(\Sym^2(V_2)) \simeq \Sym^{d-2}(E^\sharp)$.
\end{lem}

\bigskip

Lemmas~\ref{lem:base-ai} and~\ref{lem:base-qq} allow to define a linear map
$$\r_d : \Ker(\e_{2,d}) \longto \Sym^{d-2}(E^\sharp)$$
by the formula
$$\r_d(\abb_{i-1}\star\abb_{j+1}-\abb_i\star\abb_j)=
q^{\star d-j-1}\star Q^{\star i-1} \star \Psi_{j-i}^\star(\eulerb,q,Q)$$
for $1 \le i \le j \le d$. It is an isomorphism of vector spaces but 
an important fact is the following:

\bigskip

\begin{lem}\label{lem:sl2}
The map $\r_d : \Ker(\e_{2,d}) \longto \Sym^{d-2}(E^\sharp)$ is an isomorphism 
of $\Sb\Lb_2(\CM)$-modules.
\end{lem}

\bigskip

\begin{proof}
This is more or less the computation done in the end of the proof of Theorem~\ref{prop:poisson-ai-aj}. 
It is sufficient to prove that it is an isomorphism of $\sG\lG_2(\CM)$-modules. 
By~(\ref{eq:action-sl2-c}) Proposition~\ref{prop:poisson}, we have
\eqna
f \bull (\abb_{i-1}\star\abb_{j+1}-\abb_i\star\abb_j) &=& (i-1) \abb_{i-2} \star \abb_{j+1} 
+ (j+1) \abb_{i-1} \star \abb_j - i \abb_{i-1} \star \abb_j -j \abb_i \star \abb_{j-1}\\
&=& (i-1)(\abb_{i-2} \star \abb_{j+1} - \abb_{i-1} \star \abb_j) 
+ j(\abb_{i-1} \star \abb_j - \abb_i \star \abb_{j-1}).
\endeqna
Therefore, 
\eqna
\r_d(f \bull (\abb_{i-1}\star\abb_{j+1}-\abb_i\star\abb_j) &=& 
(i-1) q^{\star d-j-1}\star Q^{\star i-2} \star \Psi_{j-i+1}^\star(\eulerb,q,Q) \\
&& + j q^{\star d-j}\star Q^{\star i-1} \star \Psi_{j-i-1}^\star(\eulerb,q,Q). \\
\endeqna
and so one gets
\begin{multline*}
\r_d(f \bull (\abb_{i-1}\star\abb_{j+1}-\abb_i\star\abb_j))\\ 
= q^{\star d-j-1}\star Q^{\star i-2}  \star 
((i-1)\Psi_{j-i+1}^\star(\eulerb,q,Q) + j q\star Q \star  \Psi_{j-i-1}^\star(\eulerb,q,Q))\\
= q^{\star d-j-1}\star Q^{\star i-2}  \star ((i-1)\eulerb \star \Psi_{j-i}^\star(\eulerb,q,Q) + (j-i+1) 
q\star Q \star  \Psi_{j-i-1}^\star(\eulerb,q,Q)),
\end{multline*}
where the last equality follows from~(\ref{eq:pi-rec}). Applying now~(\ref{eq:pi-der}) yields 
\begin{multline*}
(i-1)\eulerb \star \Psi_{j-i}^\star(\eulerb,q,Q) + (j-i+1) 
q\star Q \star  \Psi_{j-i-1}^\star(\eulerb,q,Q) \\ 
= (i-1)\eulerb \star \Psi_{j-i}^\star(\eulerb,q,Q) 
+ 2 q \star Q \star \Bigl(\frac{\partial \Psi_{j-i}}{\partial T}\Bigr)^\star(\eulerb,q,Q) + 
q \star \eulerb  \Bigl(\frac{\partial \Psi_{j-i}}{\partial T'}\Bigr)^\star(\eulerb,q,Q).
\end{multline*}
Putting things together and using again~(\ref{eq:action-sl2-c}) 
and Proposition~\ref{prop:poisson} yields
$$\r_d(f \bull (\abb_{i-1}\star\abb_{j+1}-\abb_i\star\abb_j))
= f \bull \r_d(\abb_{i-1}\star\abb_{j+1}-\abb_i\star\abb_j),$$
as desired. The fact that
$$\r_d(e \bull (\abb_{i-1}\star\abb_{j+1}-\abb_i\star\abb_j))
= e \bull \r_d(\abb_{i-1}\star\abb_{j+1}-\abb_i\star\abb_j),$$
follows from a similar computation and this completes the proof of the Lemma.
\end{proof}

\bigskip

Using the isomorphism of $\Sb\Lb_2(\CM)$-modules $\r_d$, the family of 
equations~$(\ZG_{i,j})$ can be rewritten as follows:
\equat\label{eq:zij}
\forall~\ph \in \Ker(\e_{2,d}),~\ph-\r_d(\ph)\star(\eulerb^{\star 2} - 4 q \star Q - d^2a^2) 
\in \Ker(\s).
\endequat

\bigskip

\begin{rema}\label{rem:hermite}
The existence of such an isomorphism of $\Sb\Lb_2(\CM)$-modules 
$\Ker(\e_{2,d}) \longiso \Sym^{d-2}(E^\sharp)$ is a consequence 
of Hermite's reciprocity law, as it has been explained to 
us by Pierre-Louis Montagard. Indeed,
Hermite's reciprocity law (see for instance~\cite[Cor.~2.2]{brion}) 
says that we have an isomorphism of $\Sb\Lb_2(\CM)$-modules
$$h_{m,n} : \Sym^m(\Sym^n(V_2)) \longiso \Sym^n(\Sym^m(V_2))$$
making the diagram
$$\diagram
\Sym^m(\Sym^n(V_2)) \rrto^{\DS{h_{m,n}}} \drto_{\DS{\e_{m,n}}} && \Sym^n(\Sym^m(V_2)) 
\dlto^{\DS{\e_{n,m}}} \\
& \Sym^{mn}(V_2) & 
\enddiagram$$
commutative. 
In particular, $h_{m,n}$ induces an isomorphism, still denoted by $h_{m,n}$, 
between $\Ker(\e_{m,n})$ and $\Ker(\e_{n,m})$. 

In the particular case where $m=2$ and $n=d$, the kernel of the evaluation map 
$\e_{\bullet,2} : \Sym(\Sym^2(V_2))=\Sym(E^\sharp) \longto \Sym(V_2)$ is the principal ideal 
generated by $\eulerb^{\star 2} - 4 q \star Q$ so that the map 
$$\fonctio{\Sym^{d-2}(V_2)}{\Ker(\e_{d,2})}{\ph}{(\eulerb^{\star 2} - 4 q \star Q) \star \ph}$$
is an isomorphism of $\Sb\Lb_2(\CM)$-modules. Composing the inverse of this isomorphism 
with $h_{2,d}$ gives an isomorphism $\Ker(\e_{2,d}) \longiso \Sym^{d-2}(E^\sharp)$.\finl
\end{rema}

\medskip

\begin{rema}\label{rem:other-varieties}
Since $\eulerb^{\star 2} - 4 q \star Q\in\Sym(E^\sharp)^{\Sb\Lb_2(\CM)}$ (in fact, it even 
generates this invariant algebra) we can define, for any polynomial $P$ in one variable, 
a variety $\ZC^P$ by the following equations:
\begin{multline*}
\ZC^P=\{(\qG,\QG,\eG,a_0,a_1,\dots,a_d) \in \CM^{d+4}~|~\\
\forall~1\le i \le j \le d-1,~
\begin{cases}
\eG a_i = \qG a_{i+1} + \QG a_{i-1},\\
a_{i-1} a_{j+1} - a_i a_j  = 
P(\eG^2-4\qG\QG) \qG^{d-j-1} \QG^{i-1} \Psi_{j-i}(\eG,\qG,\QG)
\end{cases} \}.
\end{multline*}
By~(\ref{eq:zi-new}) and~(\ref{eq:zij}), the variety $\ZC^P$ can we rewritten as 
follows:
\begin{multline*}
\ZC^P=\{(\qG,\QG,\eG,a_0,a_1,\dots,a_d) \in \CM^{d+4}~|~\\
\begin{cases}
\forall \ph \in \Ker(\mu_{2,d})   
\cap \bigcap_{D \in \Der(\Sym(V_2))} \Ker(D^{(2)}),~\ph(\qG,\QG,\eG,a_0,a_1,\dots,a_d)=0,\\
\forall \ph \in \Ker(\e_{2,d}),~\ph(a_0,a_1,\dots,a_d)=P(\eG^2-4\qG\QG) \r_d(\ph)(\eG,\qG,\QG)
\end{cases} \}.
\end{multline*}
This shows that $\ZC^P$ is an $\Sb\Lb_2(\CM)$-stable subvariety of $\CM^{d+4} \simeq E^*$.\finl
\end{rema}

\bigskip

\section{Fixed points under diagram automorphism}\label{sec:fixe}

\medskip

Let $\sqrt{\z}$ be a primitive $2d$-th root of unity such that $(\sqrt{\z})^2=\z$ and let 
$\t=\begin{pmatrix} 0 & \sqrt{\z} \\ \sqrt{\z}^{~-1} & 0 \end{pmatrix}$. Then 
$\t s \t^{-1}=t$ and $\t t \t^{-1}=s$. So $\t$ normalizes $W$ and, since $c_s=c_t$, 
$\t$ acts on $Z_c$ and so on $\ZC_c$ by~\cite{calogero}. The action on the generators 
of $Z_c$ given in Theorem~\ref{theo:presentation-zc} is easily computed:
\equat\label{eq:action-tau}
\lexp{\t}{q}=q,\quad \lexp{\t}{Q}=Q,\quad \lexp{\t}{\eulerb}=\eulerb
\quad \text{and}\quad \lexp{\t}{\abb_i}=-\abb_i
\endequat
for $0 \le i \le d$. 

Using the description of $\ZC_c$ as a closed subvariety of $\CM^{d+4}$ 
as in~\S\ref{sub:lie} thanks to Theorem~\ref{theo:presentation-zc}, one gets:
$$\ZC_c^\t=\{(\qG,\QG,e,a_0,a_1,\dots,a_d) \in \ZC_c~|~a_0=a_1=\cdots=a_d=0\}.$$
Therefore,
$$\ZC_c^\t \simeq \{(\qG,\QG,\eG) \in \CM^3~|~\forall~1 \le i \le j \le d-1, 
(\eG^2-\qG\QG-d^2a^2)\qG^{d-j-1}\QG^{i-1}\Psi_{j-i}(e,\qG,\QG)=0\}.$$
Let $(\eG,\qG,\QG) \in \ZC_c^\t$. If $\qG \neq 0$, then the above equation with 
$i=j=1$ gives $\eG^2-\qG\QG-d^2a^2=0$. Similarly, if $\QG \neq 0$, the 
above equation with $i=j=d-1$ gives $\eG^2-\qG\QG-d^2a^2=0$. 
So assume now that $\qG=\QG=0$. Then the above equation with $i=1$ and $j=d-1$ 
gives $(\eG^2-d^2a^2)\Psi_{d-2}(e,0,0)=0$. But an easy induction on $k$ 
shows that $\Psi_k(T,0,0)=T^k$ for all $k$, so this gives 
$(\eG^2-d^2a^2)\eG^{d-2}=0$. This discussion shows that
\equat\label{eq:fix-tau}
\ZC_c^\t \simeq \{(0,0,0)\} \cup \{(\qG,\QG,\eG) \in \CM^3~|~(\eG-da)(\eG+da)=\qG\QG\}.
\endequat
So the $0$-dimensional irreducible component is of course isomorphic to 
the Calogero-Moser space associated with the trivial group (!), and the 
$2$-dimensional irreducible component is isomorphic to the Calogero-Moser 
spaces associated with the pair $(V^\t,W^\t)$ and parameter $da/2$: indeed, 
$\dim V^\t=1$, $W^\t=\langle w_0 \rangle \simeq \mub_2$ and equations 
for Calogero-Moser spaces associated with cyclic groups are given for 
instance in~\cite[Theo.~18.2.4]{calogero}. Moreover, Proposition~\ref{prop:poisson} 
shows that this isomorphism respect the Poisson bracket. So we have 
proved the following result, which confirms~\cite[Conj.~FIX]{calogero} 
(or~\cite[Conj.~B]{bonnafe auto}):

\bigskip

\begin{prop}\label{prop:fixe}
The unique $2$-dimensional irreducible component of $\ZC_c^\t$ is isomorphic, 
as a Poisson variety endowed with a $\CM^\times$-action, to the 
Calogero-Moser space associated with the pair $(V^\t,W^\t) \simeq (\CM,\mub_2)$ 
and the parameter map $\Ref(\mub_2)=\{-1\} \to \CM$, $-1 \mapsto da$.
\end{prop}

\end{document}